%&amstex

\loadeufm
\loadeusm
\loadbold
\documentstyle{amsppt}
\NoRunningHeads
\magnification=1200
\pagewidth{32pc}
\pageheight{42pc}
\vcorrection{1.2pc}

\define\wh{\widehat}

\define\bc{\Bbb C}
\define\bz{\Bbb Z}
\define\br{\Bbb R}

\define\bda{\boldkey A}
\define\bdb{\boldkey B}

\define\bdc{\boldkey C}

\topmatter
\title Hermitian representations of the extended affine Lie algebra
$\widetilde{\frak{gl}_{2}(\bc_q)}$\endtitle
\footnotetext"$^{1}$"{The author gratefully acknowledges the grant
support from the Natural Sciences and Engineering Research Council
of Canada.}
\author Yun Gao
\footnotemark"$^{1}$" \, \,  and \, \, Ziting Zeng
\endauthor
\address
Department of Mathematics and Statistics, York University,
Toronto, \newline Canada M3J 1P3
\endaddress
\email ygao\@yorku.ca \ and \ ziting\@mathstat.yorku.ca
\endemail
\dedicatory Dedicated to Professor Kyoji Saito on the occasion of
his sixtieth birthday
\enddedicatory \abstract

 We use the idea of free fields to obtain highest weight
representations for the extended affine Lie algebra
$\widetilde{\frak{gl}_{2}(\bc_q)}$ coordinatized by the quantum
torus $\bc_q$ and go on to construct a contravariant hermitian
form. We further give a necessary and sufficient condition such
that the contravariant hermitian form is positive definite.

\endabstract

\endtopmatter
\document

\subhead\S 0. Introduction\endsubhead

\medskip

Extended affine Lie algebras are a higher dimensional
generalization of affine Kac-Moody Lie algebras introduced in
[H-KT]. Even earlier than this Saito in [S] developed the notion
of extended affine root systems in the study of singularity
theory. It turns out that the non-isotropic root systems of
extended affine Lie algebras are precisely Saito's extended affine
root systems. Those Lie algebras and root systems have been
further studied in [AABGP], [BGK] and [ABGP], and among others.
There are extended affine Lie algebras which allow not only
Laurent polynomial algebra as coordinate algebra but also quantum
torus (even a nonassociative torus) depending on the type of Lie
algebra. The representations for extended affine Lie algebras and
their cousins---toroidal Lie algebras have been studied widely in
the past two decades.

In the representation theory of Lie algebras with a triangular
decomposition, the existence of a highest weight vector  and
 unitarizability are two fundamental assumptions.
Let us first recall some definitions(see [JK1]). Suppose that
$\frak g$ is a complex Lie algebra. Let $U(\frak{g})$ be its
universal enveloping algebra. Let $\Cal B$ be a subalgebra of
$\frak g$ (called a Borel subalgebra) and $\omega$ be an
antilinear anti-involution of $\frak g$ such that
$$\Cal{B} + \omega(\Cal{B}) = \frak{g}.\tag 0.1$$

Let $\lambda: \Cal{B}\to \bc$ be a 1-dimensional representation of
$\Cal{B}$. A representation $\pi: \frak{g}\to \frak{gl}(V)$ is
called a highest weight representation with highest weight
$\lambda$ if there exists a vector $v_\lambda\in V$ with the
following properties:
$$\align &\pi(U(\frak{g}))v_\lambda = V, \tag 0.2\\
&\pi(b)v_\lambda = \lambda(b) v_\lambda \text{ for any } b\in
\Cal{B}\tag 0.3\endalign $$

A hermitian form $(\cdot , \cdot)$ on $V$ such that
$$\align &(v_\lambda, v_\lambda) =1 \tag 0.4\\
&(\pi(a)u, v) = (u, \pi(\omega(a))v) \text{ for all }a\in
\frak{g}, \text{ and } u, v\in V \tag 0.5\endalign$$ is called
contravariant. One can show that, under some natural conditions,
for any highest weight $\lambda : \Cal{B}\to \bc$ there exists a
unique highest weight representation with a nondegenerate
contravariant hermitian form. As pointed out in [JK1], the
non-trivial problem is then whether this contravariant hermitian
form is positive definite(the representation $\pi$ is thus
unitarizable).

The free fields construction was first given by Wakimoto [W2] for
the affine Lie algebra $\hat{\frak{sl}}_2$ and in a great
generality by Feigin and Frenkel [FF] for the affine Lie algebras
$\hat{\frak{sl}}_n$. The book [EFK] gave a detailed treatment for
the free fields construction of the affine Lie algebra
$\hat{\frak{sl}}_2$.

In this paper we use the idea of free fields to give a new class
of highest weight representations of the extended affine Lie
algebra $\widetilde{\frak{gl}_{2}(\bc_q)}$ with respect to some
natural Borel subalgebra, where $\bc_q$ is the quantum torus (or
the algebraic version of the irrational rotation algebra in the
non-commutative geometry). This class of representations depends
on an infinite family $X$ of elements of $\text{SL}_2(\Bbb{C})$
and one complex parameter $\mu$  and is realized on the
commutative polynomial algebra $V=\Bbb{C}[x_{(m, n)} : (m,
n)\in\Bbb{Z}^2]$ in terms of the Weyl algebra
$W=\Bbb{C}[x_{(m,n)}, \frac{\partial}{\partial x_{(m,n)}}: (m,
n)\in \Bbb{Z}^2]$ twisted by an action of the family $X$ of
elements of $\text{SL}_2(\Bbb{C})$. This is the main result which
is stated in Theorem 2.12. It may be noteworthy to point out that
this realization involves operators which are cubic on standard
generators of the twisted Weyl algebra. The construction of these
representations is motivated by Wakimoto's works in particular the
unpublished manuscript [W1] where he considered the Lie algebra
$\frak{sl}_2(\bc[s^{\pm 1}, t^{\pm 1}])$. Our Theorem 3.6 provides
a contravariant hermitian form for the
$\widetilde{\frak{gl}_{2}(\bc_q)}$-module.  To find out a
necessary and sufficient condition for the contravariant hermitian
form being positive definite (see Theorem 4.8), we employ the
techniques developed by Jakobsen-Kac [JK2].

\medskip

{\it Throughout this paper, we denote the field of complex
numbers, real numbers and
 the ring of integers by $\bc$, $\br$ and $\bz$ respectively.}

\bigskip

\subhead \S 1. Extended affine Lie algebras\endsubhead

\medskip

Let $q$ be a non-zero complex number. A quantum $2$-torus associated to $q$ (see [M])
is the unital associative $\bc$-algebra $\bc_{q}[s^{\pm 1}, t^{\pm 1}]$
( or, simply $\bc_{q}$)
with generators $s^{\pm 1},  t^{\pm 1}$ and relations
$$s s^{-1} =s^{-1}s =t t^{-1}= t^{-1}t=1 \, \text{  and } \, \ ts = q st. \tag 1.1$$
 Then we have
$$(s^{m_1}t^{n_1})(s^{m_2}t^{n_2})= q^{n_1 m_2}s^{m_1+m_2}t^{n_1+n_2}\tag 1.2$$ and
$$\bc_{q} = \oplus_{m, n\in\bz} \bc s^m t^n. \tag 1.3$$

\medskip

Define $\kappa : \bc_q \to \bc$ to be  a $\bc$-linear function
given by
$$\kappa(s^m t^n) =\delta_{(m,n),(0,0)} \tag 1.4$$

Let $d_s$, $d_t$ be the degree operators on $\bc_q$ defined by
$$d_s(s^mt^n) = ms^m t^n, \, \,  d_t(s^m t^n) = n s^m t^n \tag 1.5$$
 for $m, n\in \bz$.

\medskip

For the associative algebra $\bc_q$ over $\bc$, we have the
 matrix algebra $M_{2}(\bc_q)$ with entries from $\bc_q$.
We will write $A(x)\in M_{2}(\bc_q)$ for $A\in M_{2}(\bc)$ and
$x\in \bc_q$, where $A(x)=(a_{ij}x)\in M_2(\bc_q)$  if $A=(a_{ij})
\in M_2(\Bbb{C})$. Let $\frak{gl}_{2}(\bc_q)$ be the Lie algebra
associated to $M_{2}(\bc_q)$ as usual. The Lie algebra
$\frak{gl}_{2}(\bc_q)$
 has a nondegenerate
invariant form given by
$$(A(a), B(b)) = tr(AB)\kappa(ab), \, \text{ for }A, B\in M_{2}(\bc),
a, b\in \bc_q.\tag 1.6$$

We form a  natural central extension of $\frak{gl}_{2}(\bc_q)$ as
follows.
$$\widehat{\frak{gl}_{2}(\bc_q)} = \frak{gl}_{2}(\bc_q)\oplus \bc c_s\oplus
\bc c_t \tag 1.7$$ with Lie bracket
$$\align & [A(s^{m_1}t^{n_1}), B(s^{m_2}t^{n_2})]\tag 1.8\\
=&A(s^{m_1}t^{n_1})B(s^{m_2}t^{n_2})-B(s^{m_2}t^{n_2})A(s^{m_1}t^{n_1})\\
&+tr(AB)\kappa((d_s s^{m_1}t^{n_1})s^{m_2}t^{n_2})c_s + tr(AB)
\kappa((d_t s^{m_1}t^{n_1})s^{m_2}t^{n_2})c_{t}
\endalign$$
for $m_1, m_2, n_1, n_2\in \bz, A, B\in M_{2}(\bc)$,
where $c_s$ and $c_t$ are central
elements of $\widehat{\frak{gl}_{2}(\bc_q)}$.

The derivations $d_s$ and $d_t$ can be extended
to derivations on $\frak{gl}_{2}(\bc_q)$. Now we can
define the semi-direct product of the Lie algebra $\widehat{\frak{gl}_{2}(\bc_q)}$
and those derivations:
$$\widetilde{\frak{gl}_{2}(\bc_q)} =\widehat{\frak{gl}_{2}(\bc_q)} \oplus \bc d_s
\oplus \bc d_t. \tag 1.9$$

Next we extend the nondegenerate form on $\frak{gl}_{2}(\bc_q)$ to a symmetric bilinear
form on
$\widetilde{\frak{gl}_{2}(\bc_q)}$ as follows:
$$(A(a), B(b))=tr(AB)\kappa(ab), \, (c_s, d_s) = (c_t, d_t) =1,\tag 1.10$$ all
others are zero, for $A, B\in M_{2}(\bc), a, b \in\bc_q$.

Then $\widetilde{\frak{gl}_{2}(\bc_q)}$ is an extended affine
 Lie algebra of type $A_{1}$ with nullity $2$. (See [AABGP] and
[BGK] for definitions).
% Some representations for those Lie algebras
% have been obtained in [W] for $q=1$, and [BS], [G1,2,3] for generic $q$.

Let $E_{ij}$ be the matrix whose $(i, j)$-entry is $1$ and $0$
elsewhere. Then, in $\widehat{\frak{gl}_{2}(\bc_q)}$, we have
 $$\align & [E_{ij}(s^{m_1}t^{n_1}), E_{kl}(s^{m_2}t^{n_2})]\tag 1.11\\
=&\delta_{jk}q^{n_1m_2}E_{il}(s^{m_1+m_2}t^{n_1+n_2})-\delta_{il}q^{n_2 m_1}E_{kj}
(s^{m_1+m_2}t^{n_1+n_2})\\
&+m_1q^{n_1 m_2}\delta_{jk}\delta_{il}\delta_{m_1+m_2, 0}\delta_{n_1 + n_2, 0}c_s
+n_1q^{n_1m_2}\delta_{jk}\delta_{il}\delta_{m_1+m_2, 0}\delta_{n_1+n_2, 0}c_t
\endalign$$
for $m_1, m_2, n_1, n_2\in \bz$.

The extended affine Lie algebra $\widetilde{\frak{gl}_{n}(\bc_q)}$
for $n\geq 2$ has been studied in [BGT], [BS], [E], [G1, 2, 3],
[G-KK], [VV], and among others.

\bigskip

\subhead \S 2. Representations for
 $\widetilde{\frak{gl}_{2}(\bc_q)}$
\endsubhead

\medskip

In this section, we will construct
$\widetilde{\frak{gl}_{2}(\bc_q)}$-modules by using Wakimoto's
free fields [W1, W2].  Let
$$V = \bc[x_{(m, n)}: (m,n)\in \bz^2]\tag 2.1$$
 be the (commutative) polynomial ring of infinitely many variables.
The operators $x_{(m, n)}$ and $\frac{\partial}{\partial
x_{(m,n)}}$ act on $V$ as the usual multiplication and
differentiation operators respectively.

\medskip

Given a family $X=\{X_{m, n}| (m,n) \in \Bbb{Z}^2 \}$ of $2\times
2$ lower triangular matrices, where
$$X_{m,n}=\pmatrix a_{(m, n)} & 0 \\
c_{(m, n)} & d_{(m, n)}\endpmatrix \in \text{ SL}_2(\bc) $$ for
$(m,n)\in \bz^2$ (so $a_{(m,n)} d_{(m, n)}=1$), we set
 $$\align & P_{\bda}=a_{\bda}\frac{\partial}{\partial
 x_{\bda}}\tag 2.2\\
& Q_{\bda}=c_{\bda}\frac{\partial}{\partial
x_{\bda}}+d_{\bda}x_{\bda}\tag 2.3\endalign $$ for $\bda = (m,
n)\in \bz^2$. It is easy to see the following formula holds true.
$$x_{\bda} = a_{\bda}Q_{\bda} - c_{\bda}P_{\bda}.\tag 2.4$$

\proclaim{Lemma 2.5} For $\bda, \bdb,\bdc \in \Bbb{Z}^2$, we have
$$[P_\bda,P_\bdb]=0,[Q_\bda,Q_\bdb]=0,[P_\bda,Q_\bdb]=\delta_{\bda,\bdb}.$$
\endproclaim

 \medskip

 For a fixed $\mu\in \Bbb{C}$,  define the following operators on $V$:
$$\align  e_{12}(m_1, n_1) &= -q^{-m_1n_1}\mu P_{(-m_1,-n_1)} \tag 2.6\\
&-\sum\Sb(m, n)\in\bz^2 \\
 (m', n')\in\bz^2\endSb
q^{n_1m' + nm_1 + nm'} Q_{(m+m'+m_1, n+n'+n_1)}
P_{(m,n)}P_{(m', n')}\\
e_{21}(m_1, n_1)&= Q_{(m_1, n_1)}\tag 2.7\\
e_{11}(m_1, n_1)&= -\sum_{(m, n)\in\bz^2} q^{nm_1}Q_{(m+m_1,
n+n_1)}
P_{(m,n)}-\frac{1}2\mu\delta_{(m_1,n_1),(0,0)}\tag 2.8\\
e_{22}(m_1, n_1)&= \sum_{(m, n)\in\bz^2} q^{mn_1} Q_{(m+m_1,
n+n_1)}
P_{(m,n)}+\frac{1}2\mu\delta_{(m_1,n_1),(0,0)}\tag 2.9\\
D_1 &= \sum_{(m, n)\in\bz}m Q_{(m,n)}P_{(m, n)}\tag 2.10\\
D_2 &= \sum_{(m, n)\in\bz}n Q_{(m,n)}P_{(m, n)}\tag 2.11
\endalign$$
for $m_1, n_1\in \bz$. Although $e_{11}(m_1, n_1), e_{22}(m_1,
n_1), e_{12}(m_1, n_1)$, $D_1$ and $D_2$ are infinite sums, they
are well-defined as operators on $V$.

Now we can state our first result.

 \proclaim{Theorem 2.12} The linear map $\pi_{X, \mu}:
\widetilde{\frak{gl}_{2}(\bc_q)} \to \text{End }V$ given by
$$\align & \pi_{X, \mu}(E_{ij}(s^{m_1}t^{n_1})) = e_{ij}(m_1, n_1),\\
&\pi_{X,\mu}(d_s) = D_1, \quad \pi_{X,\mu}(d_t) = D_2,\quad
\pi_{X,\mu}(c_s) = \pi_{X,\mu}(c_t) = 0,
\endalign $$
for $m_1, n_1\in \bz, 1\leq i, j\leq 2$, is a Lie algebra homomorphism.
\endproclaim
\demo{Proof} Since the parameter $q$ is involved in our
construction (2.6) through (2.9), we shall handle the
verifications in a few more details.

The following three identities are straightforward.
$$\align & [e_{11}(m_1, n_1), e_{22}(m_2, n_2)]\\
=&-\sum_{(m', n')\in\bz^2} q^{(n'+n_2)m_1 + m'n_2} Q_{(m' + m_2 + m_1, n' + n_2 + n_1)}P_{(m',n')}\\
&+\sum_{(m, n)\in\bz^2} q^{nm_1 + (m+m_1)n_2} Q_{(m + m_1 + m_2, n + n_1 + n_2)}P_{(m,n)}\\
=&0;
\endalign$$
$$\align &[e_{11}(m_1, n_1), e_{21}(m_2, n_2)]\\
=&- q^{n_2m_1}Q_{(m_1+m_2, n_1+n_2)} =-q^{n_2m_1}e_{21}(m_1+m_2,
n_1+n_2);\\
 &[e_{22}(m_1, n_1), e_{21}(m_2, n_2)]\\
 =&q^{m_2n_1}Q_{(m_1+m_2,n_1+n_2)}=q^{m_2n_1}e_{21}(m_1+m_2,
 n_1+n_2).
\endalign $$

$$\align & [e_{11}(m_1, n_1), e_{11}(m_2, n_2)]\\
=&\sum_{(m', n')\in\bz^2} q^{(n'+n_2)m_1 + n'm_2} Q_{(m' + m_2 + m_1, n' + n_2 + n_1)}P_{(m',n')}\\
&-\sum_{(m, n)\in\bz^2} q^{nm_1 + (n+n_1)m_2} Q_{(m + m_1 + m_2, n + n_1 + n_2)}P_{(m,n)}\\
=&-q^{n_2m_1}\{-\sum_{(m',
n')\in\bz^2}q^{n'(m_1+m_2)}Q_{(m'+m_1+m_2,n'+n_1+n_2)}P_{(m',n')}\\
&-\frac{1}2\mu\delta_{(m_1+m_2,n_1+n_2),(0,0)}\}\\
&+q^{n_1m_2}\{-\sum_{(m,
n)\in\bz^2}q^{n(m_1+m_2)}Q_{(m+m_1+m_2,n+n_1+n_2)}P_{(m,n)}
-\frac{1}2\mu\delta_{(m_1+m_2,n_1+n_2),(0,0)}\}\\
 =&-q^{n_2m_1}e_{11}(m_1+m_2, n_1 + n_2) +
q^{n_1m_2}e_{11}(m_1+m_2, n_1 + n_2);
\endalign$$
Similarly to the above case, one can check that
$$\align & [e_{22}(m_1, n_1), e_{22}(m_2, n_2)]\\
 =&q^{n_1m_2}e_{22}(m_1+m_2, n_1 + n_2) -
q^{n_2m_1}e_{22}(m_1+m_2, n_1 + n_2).
\endalign$$

$$\align & [e_{11}(m_1, n_1), e_{12}(m_2, n_2)]\\
=&\mu q^{-m_2n_2}\sum\Sb (m,n)\in\bz^2\endSb q^{nm_1}[Q_{(m+m_1,n+n_1)}P_{(m,n)},P_{(-m_2,-n_2)}]\\
&+\sum\Sb (m,n)\in\bz^2\\
 (m',n')\in\bz^2\\
 (m'', n'')\in\bz^2\endSb q^{nm_1 + n_2m''+ n'm_2 + n'm''}\\
 &\cdot [Q_{(m+m_1, n+n_1)}P_{(m,n)}, Q_{(m'+m''+m_2, n'+n''+n_2)}
 P_{(m',n')}P_{(m'',n'')}]\\
=&-\mu q^{-m_2n_2+(-n_1-n_2)m_1}P_{(-m_1-m_2, -n_1-n_2)}\\
+&\sum\Sb (m', n')\in\bz^2\\
 (m'', n'')\in\bz^2\endSb q^{(n'+n''+n_2)m_1 + n_2m'' + n'm_2 + n'm''}
Q_{(m' +m''+ m_2 + m_1, n'+n'' + n_2 + n_1)}P_{(m',n')}P_{(m'',n'')}\\
&-\sum\Sb (m, n)\in\bz^2\\
(m'', n'')\in\bz^2\endSb q^{nm_1 + n_2m'' + (n+n_1)m_2 +
(n+n_1)m''}
Q_{(m +m''+ m_1 + m_2, n+n''+ n_1 + n_2)}P_{(m,n)}P_{(m'',n'')}\\
&-\sum\Sb (m', n')\in\bz^2\\
 (m, n)\in\bz^2\endSb q^{nm_1 + n_2(m+m_1) + n'm_2 + n'(m+m_1)}
Q_{(m' +m+ m_1 + m_2, n'+n + n_1 + n_2)}P_{(m',n')}P_{(m,n)}\\
=&-\mu q^{-m_2n_2+(-n_1-n_2)m_1}P_{(-m_1-m_2, -n_1-n_2)}\\
&-\sum\Sb (m, n)\in\bz^2\\
(m'', n'')\in\bz^2\endSb q^{nm_1 + n_2m'' + (n+n_1)m_2 +
(n+n_1)m''}
Q_{(m +m''+ m_1 + m_2, n+n''+ n_1 + n_2)}P_{(m,n)}P_{(m'',n'')}\\
&\text{(the second and the fourth terms are negative to each other)}\\
=&q^{n_1m_2}(-\mu q^{-(m_1+m_2)(n_1+n_2)}P_{(-m_1-m_2,-n_1-n_2)}\\
&-\sum\Sb (m, n)\in\bz^2\\
(m'', n'')\in\bz^2\endSb q^{nm_1 + n_2m'' + nm_2 + (n+n_1)m''}
Q_{(m +m''+ m_1 + m_2, n+n''+ n_1 + n_2)}P_{(m,n)}P_{(m'',n'')})\\
=&q^{n_1m_2}e_{12}(m_1+m_2, n_1 + n_2);
\endalign$$
In a similar way, one may obtain that
$$[e_{22}(m_1, n_1), e_{12}(m_2, n_2)]
 =-q^{n_2m_1}e_{21}(m_1+m_2, n_1 + n_2).$$

$$\align &[e_{12}(m_1, n_1), e_{21}(m_2, n_2)]\\
=&-q^{-m_1n_1}\mu [P_{(-m_1,-n_1)},Q_{(m_2,n_2)}]\\
&-\sum\Sb (m,n)\in\bz^2\\
 (m',n')\in\bz^2\endSb q^{n_1m' + nm_1 + nm'}
[Q_{(m+m'+m_1, n+n'+n_1)}P_{(m,n)}P_{(m'n')}, Q_{(m_2, n_2)}]\\
 =&-\delta_{(-m_1,-n_1),(m_2,n_2)}q^{-m_1n_1}\mu \\
& -\sum_{(m,n)\in\bz^2} q^{n_1m_2 + nm_1 + nm_2}Q_{(m+m_2+m_1,
n+n_2+n_1)}P_{(m,n)}\\
& -\sum_{(m',n')\in\bz^2} q^{n_1m' + n_2m_1 +
n_2m'}Q_{(m'+m_2+m_1,
n'+n_2+n_1)}P_{(m'n')}\\
=&q^{n_1m_2}(-\sum_{(m,n)\in\bz^2} q^{(m_1+m_2)n}Q_{(m+m_2+m_1,
n+n_2+n_1)}P_{(m,n)}-\frac{1}2\mu \delta_{(-m_1,-n_1),(m_2,n_2)})\\
&-q^{n_2m_1}(\sum_{(m',n')\in\bz^2} q^{(n_1+n_2)m'}Q_{(m'+m_2+m_1,
n'+n_2+n_1)}P_{(m'n')}+\frac{1}2\mu
\delta_{(-m_1,-n_1),(m_2,n_2)})\\
 =& q^{n_1m_2}e_{11}(m_1+m_2, n_1+n_2) -
q^{n_2m_1}e_{22}(m_1+m_2, n_1+n_2).
\endalign $$

Next we shall handle the most complicated situation.
$$\align & [e_{12}(m_1, n_1), e_{12}(m_2, n_2)]\\
=&\sum\Sb (m,n)\in\bz^2\\
 (m',n')\in\bz^2\\
 (\bar{m}, \bar{n})\in\bz^2\\
 (\bar{m}', \bar{n}') \in\bz^2\endSb
 q^{n_1m' + nm_1 + nm'  + n_2\bar{m}' + \bar{n}m_2 + \bar{n}\bar{m}'}\\
&\cdot [Q_{(m+m'+m_1, n+n'+n_1)}P_{(m,n)}P_{(m'n')},
Q_{(\bar{m}+\bar{m}'+m_2, \bar{n}+\bar{n}'+n_2)}P_{(\bar{m},\bar{n})}P_{(\bar{m'},\bar{n'})}]\\
&+q^{(-m_1n_1)}\mu\sum\Sb (m,n)\in\bz^2\\
 (m',n')\in\bz^2\endSb q^{n_2m'+nm_2+nm'}[P_{(-m_1,-n_1)},Q_{(m+m'+m_2,n+n'+n_2)}P_{(m,n)}P_{(m',n')}]\\
&+q^{(-m_2n_2)}\mu\sum\Sb (m,n)\in\bz^2\\
 (m',n')\in\bz^2\endSb q^{n_2m'+nm_2+nm'}[Q_{(m+m'+m_1,n+n'+n_1)}P_{(m,n)}P_{(m',n')},P_{(-m_2,-n_2)}]\\
 =&J_1 + J_2 - J_3 - J_4+J_5+J_6
\endalign $$
where
$$\align J_1 &= \sum\Sb (m,n)\in\bz^2\\
 (\bar{m}, \bar{n})\in\bz^2\\
 (\bar{m}', \bar{n}') \in\bz^2\endSb
 q^{n_1(\bar{m}+\bar{m}'+m_2) + nm_1 + n(\bar{m}+\bar{m}'+m_2)  + n_2\bar{m}' + \bar{n}m_2
+ \bar{n}\bar{m}'}\\
&\cdot
Q_{(m+\bar{m}+\bar{m}'+m_2+m_1, n+\bar{n}+\bar{n}'+n_2+n_1)}P_{(m,n)}P_{(\bar{m},\bar{n})}P_{(\bar{m}',\bar{n}')}\\
 J_2 &= \sum\Sb (m',n')\in\bz^2\\
 (\bar{m}, \bar{n})\in\bz^2\\
 (\bar{m}', \bar{n}') \in\bz^2\endSb
 q^{n_1m'+ (\bar{n}+\bar{n}'+n_2)m_1 + (\bar{n}+\bar{n}'+n_2)m'  + n_2\bar{m}' + \bar{n}m_2
+ \bar{n}\bar{m}'}\\
&\cdot
Q_{(\bar{m}+\bar{m}'+m'+m_2+m_1, \bar{n}+\bar{n}'+n'+n_2+n_1)}P_{(m',n')}P_{(\bar{m},\bar{n})}P_{(\bar{m}',\bar{n}')}\\
J_3 &= \sum\Sb (m,n)\in\bz^2\\
 (\bar{m}, \bar{n})\in\bz^2\\
 ({m}', {n}') \in\bz^2\endSb
 q^{n_1{m}' + nm_1
+ n{m}' + n_2(m+{m}'+m_1) + \bar{n}m_2 + \bar{n}(m+{m}'+m_1)}\\
&\cdot Q_{(\bar{m}+m+{m}'+m_2+m_1, \bar{n}+ n + {n}'+n_2+n_1)}P_{(\bar{m},\bar{n})}P_{(m,n)}P_{(m',n')}\\
J_4 &= \sum\Sb (m',n')\in\bz^2\\
 (m, n)\in\bz^2\\
 (\bar{m}', \bar{n}') \in\bz^2\endSb
 q^{n_1{m}' + nm_1
+ n {m}' + n_2\bar{m}'+ (n+{n}'+n_1)m_2 + (n+ {n}'+n_1)\bar{m}'}\\
&\cdot Q_{(m+m'+\bar{m}'+m_2+m_1,
n+{n}'+\bar{n}'+n_2+n_1)}P_{(\bar{m}',\bar{n}')}P_{(m,n)}P_{(m',n')}\\
J_5 &=q^{-m_1n_1}\mu\sum\Sb (m,n)\in\bz^2\endSb
q^{n_2(-m_1-m_2-m+nm_2+n(-m_1-m_2-m))}\\
&\cdot P_{(-m_1-m_2-m,-n_1-n_2-n)}P_{(m,n)}\\
J_6 &=-q^{-m_2n_2}\mu\sum\Sb (m',n')\in\bz^2\endSb
q^{n_1m'+(-n_1-n_2-n')m_1+(-n_1-n_2-n')m'}\\
&\cdot P_{(-m_1-m_2-m',-n_1-n_2-n')}P_{(m',n')}
\endalign$$
Note that $J_1 = J_4$, $J_2=J_3$ and $J_5=-J_6$. Thus
$$[e_{12}(m_1, n_1), e_{12}(m_2, n_2)]= 0.$$

It is clear that $[e_{21}(m_1, n_1), e_{21}(m_2, n_2)]= 0$. Next
we check the identities involving $D_1$ and $ D_2$.

It is obvious that the following identities hold.
$$[D_1, D_2]=0,\
[D_1, e_{21}(m_1, n_1)] = m_1e_{21}(m_1, n_1).$$

$$\align & [D_1, e_{12}(m_1, n_1)]\\
=&-q^{-m_1n_1}\mu\sum\Sb(m,n)\in\bz^2\endSb m
[Q_{(m,n)}P_{(m,n)},P_{(-m_1,-n_1)}]\\
&-\sum\Sb (m,n)\in\bz^2\\
 (m',n')\in\bz^2\\
 (m'', n'')\in\bz^2\endSb mq^{n_1m'' + n'm_1 + n'm''}
[Q_{(m, n)}P_{(m,n)}, Q_{(m'+m''+m_1, n'+n''+n_1)}P_{(m', n')}P_{(m'',n'')}]\\
=&q^{-m_1n_1}\mu (-m_1)P_{(-m_1,-n_1)}\\
&-\sum\Sb (m', n')\in\bz^2\\
 (m'', n'')\in\bz^2\endSb (m'+m''+m_1)q^{n_1m'' + n'm_1 + n'm''}
Q_{(m' +m''+ m_1, n'+n'' + n_1)}P_{(m',n')}P_{(m'',n'')}\\
&+\sum\Sb (m, n)\in\bz^2\\
(m'', n'')\in\bz^2\endSb mq^{n_1m'' + nm_1 + nm''}
Q_{(m +m''+ m_1, n+n'' + n_1)}P_{(m,n)}P_{(m'',n'')}\\
&+\sum\Sb (m', n')\in\bz^2\\
 (m, n)\in\bz^2\endSb mq^{ n_1m + n'm_1 + n'm}
Q_{(m' +m+ m_1, n'+n + n_1)}P_{(m',n')}P_{(m,n)}\\
=&-m_1q^{-m_1n_1}\mu P_{(-m_1,-n_1)}\\
& -m_1\sum\Sb (m', n')\in\bz^2\\
 (m, n)\in\bz^2\endSb q^{ n_1m + n'm_1 + n'm}
Q_{(m' +m+ m_1, n'+n + n_1)}P_{(m',n')}P_{(m,n)}\\
 =& m_1 e_{12}(m_1, n_1).
\endalign$$

$$\align &[D_1, e_{11}(m_1, n_1)]\\
=&-\sum\Sb(m,n)\in\bz^2 \\
 (m',n')\in\bz^2\endSb mq^{n'm_1}[Q_{(m, n)}P_{(m,n)}, Q_{(m'+m_1, n'+n_1)}P_{(m',n')}]\\
=&-\sum_{(m', n')\in\bz^2} (m'+m_1)q^{n'm_1} Q_{(m' + m_1, n' + n_1)}P_{(m',n')}\\
&+\sum_{(m, n)\in\bz^2}m q^{nm_1} Q_{(m + m_1, n + n_1)}P_{(m,n)}\\
=& m_1(-\sum_{(m, n)\in\bz^2}m q^{nm_1} Q_{(m + m_1, n +
n_1)}P_{(m,n)}-\frac{1}2\mu\delta_{(m_1,n_1),(0,0)})=m_1e_{11}(m_1,
n_1);
\endalign$$
Similarly to the above case, one has
$$[D_1, e_{22}(m_1, n_1)]
 =m_1e_{22}(m_1, n_1).$$

Replacing $D_1$ by $D_2$ in the above proof,  one can show that
$$\align &[D_2, e_{21}(m_1, n_1)] = n_1 e_{21}(m_1, n_1), \ [D_2, e_{11}(m_1, n_1)]= n_1 e_{11}(m_1, n_1),\\
&[D_2, e_{22}(m_1, n_1)] = n_1e_{22}(m_1, n_1), \  [D_2,
e_{12}(m_1, n_1)] = n_1 e_{12}(m_1, n_1).
\endalign$$

Therefore by comparing with (1.11) we see that the linear map
$\pi_{X, \mu}$ is indeed a Lie algebra homomorphism.  \qed\enddemo

\bigskip

\subhead \S 3. Hermitian forms \endsubhead

\medskip

Here we shall unify the hermitian forms independently studied by
Wakimoto [W1] and Jakobsen-Kac [JK2].

Define $\omega:\widetilde{\frak{gl}_{2}(\bc_q)}\mapsto
\widetilde{\frak{gl}_{2}(\bc_q)}$ a $\br$-linear map as the
following:
$$\align \omega(\lambda x)&=\bar{\lambda}\omega (x),\forall
\lambda\in \bc,x\in \widetilde{\frak{gl}_{2}(\bc_q)} \tag 3.1 \\
\omega(E_{ij}(a))&=(-1)^{i+j}E_{ji}(\overline{a}), a \in \bc_q \tag 3.2\\
\omega(d_s)&=d_s, \, \omega(d_t)=d_t, \, \omega(c_s)=c_s, \,
\omega(c_t)=c_t \tag 3.3
   \endalign$$
where $\br-$linear map $\bar{ }: \bc_q\rightarrow \bc_q$ is
defined as $\overline{\lambda
s^mt^n}=\bar{\lambda}t^{-n}s^{-m}=\bar{\lambda}q^{mn}s^{-m}t^{-n}$,
and $\bar{\lambda}$ is the complex conjugate, for any $\lambda\in
\bc$, and $m,n\in \bz$.

In the following sections, we always assume that $q\bar{q}=1$ (or
$|q| = 1$). This assumption will guarantee that the map \ $\bar{
}$ \ is of order two.

 \proclaim{Lemma 3.4}$\omega$ is an anti-linear anti-involution of
$\widetilde{\frak{gl}_{2}(\bc_q)}$.
\endproclaim
\demo{Proof} Since
$$\omega(E_{ij}(s^mt^n))=(-1)^{i+j}q^{mn}E_{ji}(s^{-m}t^{-n}), \tag 3.5 $$ we
have
$$\align
&\omega^2(E_{ij}(s^mt^n))=\omega((-1)^{i+j}q^{mn}E_{ji}(s^{-m}t^{-n}))\\
=&(-1)^{i+j}\bar{q}^{mn}(-1)^{i+j}q^{mn}E_{ij}(s^mt^n)=E_{ij}(s^mt^n),\endalign$$
so $\omega^2=id$. We only need to check
$\omega([a,b])=[\omega(b),\omega(a)]$, for any $a,b\in
\widetilde{\frak{gl}_{2}(\bc_q)}$.

$$\align &[\omega(E_{ij}(s^{m_1}t^{n_1})),\omega(E_{kl}(s^{m_2}t^{n_2}))]\\
=&(-1)^{j+k}q^{m_1n_1+m_2n_2+m_2n_1}\delta_{il}E_{jk}(s^{-(m_1+m_2)}t^{-(n_1+n_2)})\\
&-(-1)^{l+i}q^{m_1n_1+m_2n_2+m_1n_2}\delta_{jk}E_{li}
(s^{-(m_1+m_2)}t^{-(n_1+n_2)})\\
&-q^{m_1n_1+m_2n_2+m_2n_1}m_1\delta_{jk}\delta_{il}\delta_{n_1+n_2,
0}\delta_{m_1 + m_2,
0}c_s\\
&-q^{m_1n_1+m_2n_2+m_2n_1}n_1\delta_{jk}\delta_{il}\delta_{n_1+n_2,
0}\delta_{m_1+m_2, 0}c_t
\endalign$$
Thus
$$\align &\omega([E_{ij}(s^{m_1}t^{n_1}),E_{kl}(s^{m_2}t^{n_2})])\\
=&\bar{q}^{m_2n_1}\delta_{jk}\omega(E_{il}(s^{m_1+m_2}t^{n_1+n_2}))-\bar{q}^{
m_1n_2}\delta_{il}\omega(E_{kj}
(s^{m_1+m_2}t^{n_1+n_2}))\\
&+m_1\bar{q}^{m_2n_1}\delta_{jk}\delta_{il}\delta_{m_1+m_2,
0}\delta_{n_1 + n_2, 0}\omega(c_s)
+n_1\bar{q}^{m_2n_1}\delta_{jk}\delta_{il}\delta_{m_1+m_2,
0}\delta_{n_1+n_2, 0}\omega(c_t)\\
=&(-1)^{l+i}q^{m_1n_1+m_2n_2+m_1n_2}\delta_{jk}E_{li}(s^{-(m_1+m_2)}t^{-(n_1+n_2})\\
&-(-1)^{j+k}q^{m_1n_1+m_2n_2+m_2n_1}\delta_{il}E_{jk}
(s^{-(m_1+m_2)}t^{-(n_1+n_2)})\\
&+m_1q^{m_1n_1+m_2n_2+m_2n_1}\delta_{jk}\delta_{il}\delta_{m_1+m_2,
0}\delta_{n_1 + n_2, 0}c_s\\
&+n_1q^{m_1n_1+m_2n_2+m_2n_1}\delta_{jk}\delta_{il}\delta_{m_1+m_2,
0}\delta_{n_1+n_2, 0}c_t\\
=&-[\omega(E_{ij}(s^{m_1}t^{n_1})),\omega(E_{kl}(s^{m_2}t^{n_2}))]=[\omega(E_{kl}(s^{m_2}t^{n_2})),\omega(E_{ij}(s^{m_1}t^{n_1}))].
\endalign$$

As for identities involving $d_s$ and $d_t$, we have
$$[\omega(d_s),\omega(E_{ij}(s^mt^n))]=[d_s,(-1)^{j+i}q^{mn}E_{ji}(s^{-m}t^{-n})]=-(-1)^{j+i}mq^{mn}E_{ji}(s^{-m}t^{-n})$$
$$\omega([d_s,E_{ij}(s^mt^n)])=\omega(mE_{ij}(s^mt^n))=(-1)^{j+i}mq^{mn}E_{ji}(s^{-m}t^{-n})$$
Hence we get
$[\omega(d_s),\omega(E_{ij}(s^mt^n))]=\omega(-[d_s,E_{ij}(s^mt^n)])$.
Similarly,
$$[\omega(d_t),\omega(E_{ij}(s^mt^n))]=\omega(-[d_t,E_{ij}(s^mt^n)])$$
The other cases are  trivial and so the proof is completed.
\qed\enddemo

\proclaim{Theorem 3.6} Assume that $\mu$ is a real number. Then
there exists a contravariant, with respect to $\pi_{X,\mu}$ and
$\omega$,  hermitian form $(\cdot , \cdot)$ on $V$ so that
$$(\pi_{X, \mu} (a).f,g)=(f,\pi_{X, \mu}(\omega(a)).g)\tag 3.7 $$
for every $f,g\in V,a\in \widetilde{\frak{gl}_{2}(\bc_q)}$.
\endproclaim
\demo{Proof} Since $V$ is a polynomial algebra, it is sufficient
to define the form on a pair of monomials in the variables
$x_{(m,n)}, (m, n)\in \Bbb{Z}^2$. Given a monomial
$$f=\prod\Sb (m,n)\in \bz^2\endSb x_{(m,n)}^{\bda_{(m,n)}}
\tag 3.8$$ whose degree in $x_{(m, n)}$ is positive, where
$\bda_{(m,n)}\in \bz_+\bigcup \{0\}$ and only finitely many
$\bda_{(m,n)}\neq 0$, denote by $f_{\widehat{(m,n)}}$ the unique
monomial such that
$$f =x_{(m,n)}f_{\widehat{(m,n)}}. \tag 3.9$$
Denote by deg $f$ the total degree of $f$.

Now we define a hermitian form $(f,g)$ on $V$ inductively on the
degree of $f$. Since a hermitian form requires
$(f,g)=\overline{(g,f)}$, we only need to define $(f,g)$ with
$deg(g)\leqq deg(f)$.

We set $$\align & (1, 1) =1, \tag 3.10\\
 & (x_{(m,n)},1)=0, \tag 3.11\\
& (x_{(m,n)},x_{(l,k)})=\mu a_{(m,n)}\overline{a_{(m,n)}}
\delta_{(m,n),(l,k)}.\tag 3.12 \endalign $$

 Fix a positive integer $N$ and assume that the form is defined for all monomials $f$, $g$
such that $deg(g),  deg(f)\leqq N-1$ and satisfies
$$(\pi_{X, \mu}(E_{21}(s^mt^n)).f,g)=(f,\pi_{X, \mu}(\omega(E_{21}(s^mt^n))).g)\tag 3.13$$
with $deg(g)\leqq N-2, deg(f)\leqq N-1$. It is easy to see that
(3.13) holds true when deg $f$, deg $g$ $\leq 1$. Take $f$ with
$deg(f)=N$, and choose $(m,n)\in\bz^2$ such that the degree of $f$
in $x_{(m,n)} \geq 1$.

Observe that
$$f= x_{(m,n)}f_{\widehat{(m,n)}}=a_{(m,n)}Q_{(m,n)}f_{\widehat{(m,n)}}-c_{(m,n)}P_{(m,n)}f_{\widehat{(m,n)}}$$
by (2.4). Since $Q_{(m,n)}=\pi_{X,\mu}(E_{21}(s^mt^n))$ this can
be written as
$$f =a_{(m,n)}\pi_{X,\mu}(E_{21}(s^mt^n))f_{\widehat{(m,n)}}-c_{(m,n)}P_{(m,n)}f_{\widehat{(m,n)}}\tag 3.14$$

Suppose first that deg $g < N$. Then set
$$(f, g):= a_{(m,n)}(f_{\widehat{(m,n)}},
\pi_{X,\mu}(\omega(E_{21}(s^mt^n))g) -
c_{(m,n)}(P_{(m,n)}f_{\widehat{(m,n)}}, g).\tag 3.15$$ Note that
all terms here are defined by induction and (3.13) holds.

Suppose now that deg $g=N$. The first term in (3.15) still makes
sense as deg $\pi_{X, \mu}(\omega(E_{21}(s^mt^n))g <$ deg $g$. On
the other hand, $(g, P_{(m,n)}f_{\widehat{(m,n)}})$ is defined by
(3.15) as deg $P_{(m,n)}f_{\widehat{(m,n)}} <$ deg
$f_{\widehat{(m,n)}} < N$. Then the last term is also defined by
applying the same formula to $g$ and
$P_{(m,n)}f_{\widehat{(m,n)}}$ and using the fact that the form is
hermitian.

\medskip

 We have to show (3.15) is well-defined, which
means that the right-hand side of (3.15) is independent of the
choice of $(m,n)$. Namely,  we need to show that if
$\bda_{(m,n)}\geq 1$, $\bda_{(l,k)}\geq 1$ and $(m,n) \neq (l,
k)$, we have
$$\align &a_{(m,n)}(f_{\widehat{(m,n)}}, -q^{mn}e_{12}(-m,-n).g)-
c_{(m,n)}(P_{(m,n)}f_{\widehat{(m,n)}}, g)\tag 3.16\\
=&a_{(l,k)}(f_{\widehat{(l,k)}},-q^{lk}e_{12}(-l,-k).g)-c_{(l,k)}(P_{(l,k)}f_{\widehat{(l,k)}},g).
\endalign$$

Since
$$f_{\widehat{(m,n)}}=a_{(l,k)}e_{21}(l,k).f_{\widehat{(m,n)}\widehat{(l,k)}}-
c_{(l,k)}P_{(l,k)}f_{\widehat{(m,n)}\widehat{(l,k)}},\tag 3.17$$
substituting to the left-hand side of (3.16), we obtain
$$\align & \text{ LHS of (3.16)}\tag 3.18\\
=&a_{(m,n)}(a_{(l,k)}e_{21}(l,k).f_{\widehat{(m,n)}\widehat{(l,k)}}-c_{(l,k)}P_{(l,k)}f_{\widehat{(m,n)}
\widehat{(l,k)}}, -q^{mn}e_{12}(-m,-n).g)\\
&-c_{(m,n)}(P_{(m,n)}(a_{(l,k)}e_{21}(l,k).f_{\widehat{(m,n)}\widehat{(l,k)}}-c_{(l,k)}
P_{(l,k)}f_{\widehat{(m,n)}\widehat{(l,k)}}),g)\\
=&a_{(m,n)}
a_{(l,k)}(e_{21}(l,k).f_{\widehat{(m,n)}\widehat{(l,k)}},-q^{mn}e_{12}(-m,-n).g)\\
&-a_{(m,n)}c_{(l,k)}(P_{(l,k)}f_{\widehat{(m,n)}\widehat{(l,k)}},-q^{mn}e_{12}(-m,-n).g)\\
&-c_{(m,n)}a_{(l,k)}(P_{(m,n)}e_{21}(l,k).f_{\widehat{(m,n)}\widehat{(l,k)}},g)\\
&+c_{(m,n)}c_{(l,k)}(P_{(m,n)}P_{(l,k)}f_{\widehat{(m,n)}\widehat{(l,k)}},g)\\
=&a_{(m,n)}
a_{(l,k)}(f_{\widehat{(m,n)}\widehat{(l,k)}},(-1)q^{lk}e_{12}(-l,-k).(-1)q^{mn}e_{12}(-m,-n).g)\\
&-a_{(m,n)}c_{(l,k)}(e_{21}(m,n)P_{(l,k)}f_{\widehat{(m,n)}\widehat{(l,k)}},g)\\
&-c_{(m,n)}a_{(l,k)}(P_{(m,n)}e_{21}(l,k).f_{\widehat{(m,n)}\widehat{(l,k)}},g)\\
&+c_{(m,n)}c_{(l,k)}(P_{(m,n)}P_{(l,k)}f_{\widehat{(m,n)}\widehat{(l,k)}},g)
\endalign$$

Exchanging $(m,n)$ and $(l,k)$ in $(3.18)$ and noting that
$f_{\widehat{(m,n)}\widehat{(l,k)}}=f_{\widehat{(l,k)}\widehat{(m,n)}}$,
we get the right-hand side of $(3.16)$:
$$\align & \text{RHS of (3.16)}\tag 3.19\\
=&a_{(l,k)}
a_{(m,n)}(f_{\widehat{(m,n)}\widehat{(l,k)}},q^{mn}e_{12}(-m,-n).q^{lk}e_{12}(-p,-q).g)\\
&-a_{(l,k)}c_{(m,n)}(e_{21}(l,k)P_{(m,n)}f_{\widehat{(m,n)}\widehat{(l,k)}},g)\\
&-c_{(l,k)}a_{(m,n)}(P_{(l,k)}e_{21}(m,n).f_{\widehat{(m,n)}\widehat{(l,k)}},g)\\
&+c_{(l,k)}c_{(m,n)}(P_{(l,k)}P_{(m,n)}f_{\widehat{(m,n)}\widehat{(l,k)}},g)
\endalign$$

Since $[e_{12}(-m,-n),e_{12}(-l,-k)]=[P_{(m,n)},P_{(l,k)}]=0$,
subtracting (3.18) from (3.19) we have
$$\align (3.19)-(3.18)\\
=&-a_{(l,k)}c_{(m,n)}([Q_{(l,k)},P_{(m,n)}].f_{\widehat{(m,n)}\widehat{(l,k)}},g)\\
&-c_{(l,k)}a_{(m,n)}([P_{(l,k)},Q_{(m,n)}].f_{\widehat{(m,n)}\widehat{(l,k)}},g)\\
=&\delta_{(l,k),(m,n)}(a_{(l,k)}c_{(m,n)}-c_{(l,k)}a_{(m,n)})(f_{\widehat{(m,n)}\widehat{(l,k)}},g)\\
=&0
\endalign$$

Hence (3.16) holds true and (3.15) is well-defined. So we obtained
a form on $V$, and the form satisfies (3.13) for any $f,g \in V$.

Since (3.13) holds and $E_{22}(x)$ is a linear combination of
$E_{11}(x)$ and $[E_{12}(x'), E_{21}(x'')]$, in order to prove
that the form we defined is contravariant it remains to check
$$\align &(\pi_{X, \mu}(d_s).f,g)=(f,\pi_{X, \mu}(\omega(d_s)).g), \tag 3.20\\
&(\pi_{X, \mu}(d_t).f,g)=(f,\pi_{X, \mu}(\omega(d_t)).g), \tag 3.21 \\
&(\pi_{X, \mu}(E_{11}(s^mt^n)).f,g)=(f,\pi_{X,
\mu}(\omega(E_{11}(s^mt^n))).g). \tag 3.22 \endalign $$

 We do this by using induction on the
degree of $f$ and $g$.

 First we have
 $$\align &(\pi_{X,\mu}(d_s).1,1)=(1,\pi_{X, \mu}(\omega(d_s)).1)=0,\\
 &(\pi_{X, \mu}(d_s).x_{(m,n)},1)=0=(x_{(m,n)},\pi_{X, \mu}(d_s).1), \\
&(\pi_{X,
\mu}(d_s).x_{(m,n)},x_{(l,k)})=m(x_{(m,n)},x_{(l,k)})=l(x_{(m,n)},x_{(l,k)})\\
&=(x_{(m,n)},\pi_{X, \mu}(\omega(d_s)).x_{(l,k)})\endalign$$ for
$(m, n), (l, k)\in\bz^2$.

Assuming
 for any $deg(g)\leqq deg(f)\leqq N-1$, we have
 $$(\pi_{X, \mu}(d_s).f,g)=(f,\pi_{X, \mu}(\omega(d_s)).g). $$
 According to (3.14) or (2.4),
$$f=a_{(m,n)}e_{21}(m,n).h-c_{(m,n)}P_{(m,n)}h \tag 3.23$$
where $ h=f_{\widehat{(m,n)}} $  and $deg(h)=N-1,$ together with
the assumption, then
$$\align &(\pi_{X,\mu}(d_s).f,g)\\
=&a_{(m,n)}(D_1e_{21}(m,n)).h,g)-c_{(m,n)}(D_1P_{(m,n)}h,g)\\
=&a_{(m,n)}(e_{21}(m,n)D_1.h+[D_1,e_{21}(m,n)].h,g)-c_{(m,n)}(P_{(m,n)}h,D_1.g)\\
=&a_{(m,n)}(D_1.h,-q^{mn}e_{12}(-m,-n).g)+a_{(m,n)}(me_{21}(m,n).h,g)\\
&-c_{(m,n)}(P_{(m,n)}.h,D_1.g)\\
=&a_{(m,n)}(h,-q^{mn}D_1e_{12}(-m,-n).g)+a_{(m,n)}m(h,-q^{(mn)}e_{12}(-m,-n).g)\\
&-c_{(m,n)}(P_{(m,n)}.h,D_1.g)\\
=&a_{(m,n)}(h,-q^{mn}(e_{12}(-m,-n)D_1+[D_1,e_{12}(-m,-n)]).g)\\
&+a_{(m,n)}m(h,-q^{(mn)}e_{12}(-m,-n).g)-c_{(m,n)}(P_{(m,n)}.h,D_1.g)\\
=&a_{(m,n)}(h,-q^{mn}e_{12}(-m,-n)D_1.g)+a_{(m,n)}(h,+mq^{mn}e_{12}(-m,-n).g)\\
&+a_{(m,n)}m(h,-q^{(mn)}e_{12}(-m,-n).g)-c_{(m,n)}(P_{(m,n)}.h,D_1.g)\\
=&a_{(m,n)}(e_{21}(m,n).h,D_1.g)-c_{(m,n)}(P_{(m,n)}.h,D_1.g)\\
=&(f,D_1.g) = (f,\pi_{X,\mu}(d_s).g).
\endalign$$
Therefore (3.20) holds true and so does (3.21).

As for (3.22) we first have
$$\align
&(\pi_{X, \mu}(E_{11}(s^lt^k)).1,1)=-\frac{1}{2}\mu\delta_{(l,k),(0,0)}=(1,\pi_{X, \mu}(\omega(E_{11}(s^lt^k))).1),\\
&(\pi_{X, \mu}(E_{11}(s^lt^k)).x_{m,n},1)=0=(x_{m,n},\pi_{X, \mu}(\omega(E_{11}(s^lt^k))).1).\\
\endalign$$

 Secondly,
$$\align&(\pi_{X, \mu}(E_{11}(s^lt^k)).x_{(m_1,n_1)},x_{(m_2,n_2)})\\
=&-q^{ln_1}a_{(m_1,n_1)}d_{(m_1+l,n_1+k)}(x_{(m_1+l,n_1+k)},x_{(m_2,n_2)})-\frac{1}{2}\mu
\delta_{(l,k),(0,0)}(x_{(m_1,n_1)},x_{(m_2,n_2)})\\
 =&-q^{ln_1}a_{(m_1,n_1)}
\overline{a_{(m_2,n_2)}}\mu\delta_{(m_1+l,n_1+k),(m_2,n_2)}-\frac{1}{2}|a_{(m_2,n_2)}|^2
\mu^2\delta_{(l,k),(0,0)}\delta_{(m_1,n_1),(m_2,n_2)}\endalign$$
and
$$\align&(x_{(m_1,n_1)},\pi_{X, \mu}(\omega(E_{11}(s^lt^k))).x_{(m_2,n_2)})\\
=&-\overline{q^{lk-ln_2}}\overline{a_{(m_2,n_2)}}\overline{d_{(m_2-l,n_2-k)}}(x_{(m_1,n_1)},x_{(m_2-l,n_2-k)})\\
&-\overline{q^{lk}}\frac{1}{2}\mu
\delta_{(-l,-k),(0,0)}(x_{(m_1,n_1)},x_{(m_2,n_2)})\\
=&-q^{ln_1}a_{(m_1,n_1)}
\overline{a_{(m_2,n_2)}}\mu\delta_{(m_1+l,n_1+k),(m_2,n_2)}-\frac{1}{2}|a_{(m_2,n_2)}|^2
\mu^2\delta_{(l,k),(0,0)}\delta_{(m_1,n_1),(m_2,n_2)}
\endalign$$ yields
$$(\pi_{X, \mu}(E_{11}(s^lt^k)).x_{(m_1,n_1)},x_{(m_2,n_2)})=(x_{(m_1,n_1)},\pi_{X, \mu}(\omega(E_{11}(s^lt^k))).x_{(m_2,n_2)}).$$

Now assume for any $deg(g)\leqq deg(f)\leqq N-1$, we have
$$(\pi_{X, \mu}(E_{11}(s^mt^n)).f,g)=(f,\pi_{X, \mu}(\omega(E_{11}(s^mt^n))).g).$$
According to $(3.23)$, we have
$$\align &(\pi_{X,\mu}(E_{11}(s^lt^k)).f,g)\\
=&a_{(m,n)}(e_{11}(l,k)e_{21}(m,n).h,g)
-c_{(m,n)}(e_{11}(l,k)P_{(m,n)}h,g)\\
=&a_{(m,n)}(e_{21}(m,n)e_{11}(l,k).h,g)+a_{(m,n)}([e_{11}(l,k),e_{21}(m,n)].h,g)\\
&-c_{(m,n)}(P_{(m,n)}.h,q^{lk}e_{11}(-l,-k).g)\\
=&a_{(m,n)}(e_{11}(l,k).h,q^{mn}e_{12}(-m,-n)g)+a_{(m,n)}(-q^{nl}e_{21}(m+l,n+k).h,g)\\
&-c_{(m,n)}(P_{(m,n)}.h,q^{lk}e_{11}(-l,-k).g)\\
=&a_{(m,n)}(h,q^{mn}q^{lk}e_{11}(-l,-k)e_{12}(-m,-n).g)\\
&+a_{(m,n)}(h,-\bar{q}^{nl}q^{(m+l)(n+k)}e_{12}(-m-l,-n-k).g)\\
&-c_{(m,n)}(P_{(m,n)}.h,q^{lk}e_{11}(-l,-k).g)\\
=&a_{(m,n)}(h,q^{mn+lk}e_{12}(-m,-n)e_{11}(-l,-k).g)\\
&+a_{(m,n)}(h,q^{mn+lk}[e_{11}(-l,-k),e_{12}(-m,-n)].g)\\
&+a_{(m,n)}(h,-q^{-nl}q^{(m+l)(n+k)}e_{12}(-m-l,-n-k).g)\\
&-c_{(m,n)}(P_{(m,n)}.h,q^{lk}e_{11}(-l,-k)g) \\
=&a_{(m,n)}(e_{21}(m,n).h,q^{lk}e_{11}(-l,-k).g)+a_{(m,n)}(h,q^{mn+lk}q^{km}e_{12}(-m-l,-n-k).g)\\
&+a_{(m,n)}(h,-q^{mn+lk+mk}e_{12}(-m-l,-n-k).g)-c_{(m,n)}(P_{(m,n)}h,q^{lk}e_{11}(-l,-k)g) \\
=&a_{(m,n)}(e_{21}(m,n).h,q^{lk}e_{11}(-l,-k).g)-c_{(m,n)}(P_{(m,n)}.h,q^{lk}e_{11}(-l,-k)g) \\
=&(f,q^{lk}e_{11}(-l,-k)g)=(f,\pi_{X,\mu}(\omega(E_{11}(s^lt^k))).g)
\endalign$$

Hence (3.22) is also true and the form is indeed a contravariant
hermitian form on $V$. \qed
 \enddemo

\bigskip

\subhead \S 4. Conditions  for unitarity \endsubhead

\medskip

It is important to have the contravariant hermitian form on $V$ to
be positive definite so that the underlying module  is
unitarizable.

>From the definition of our contravariant form in Theorem 3.6, we
see that
$$\align(x_{(m,n)}^2, 1)=&a_{(m,n)}(x_{(m,n)},q^{mn}e_{12}(-m,-n).1)-c_{(m,n)}(P_{(m,n)}.x_{(m,n)},1)\\
=&-a_{(m,n )}c_{(m,n)}.\endalign$$ Thus the hermitian form on
different degrees can be non-zero.  It is therefore difficult to
determine when the hermitian form is positive definite. For this
we use another base of $V$ as in [JK2] rather than the natural
monomial base of $V$. We further work out the necessary and
sufficient conditions for the unitarity. We shall follow the
approach in [JK2].

\proclaim{Definition 4.1} If the hermitian form is positive
definite, $\pi_{X, \mu}$ is said to be unitarizable (w.r.t
$\omega$).
\endproclaim

Here we simplify $\pi_{X, \mu}(E_{ij}(r)).v$ as $E_{ij}(r).v$, for
any $v\in V, r\in \bc_q$.

\proclaim{Lemma 4.2} The elements
$E_{21}(r_1)E_{21}(r_2)...E_{21}(r_k).1$, where $k\in
\bz_+\bigcup\{0\}$, $r_i=s^{m_i}t^{n_i}$, $i=1...n,m_i,n_i\in \bz$
forms a basis for $V$. Moreover, if $f$ is a monomial of degree
$N$, then $f$ can be written as a linear combination of
$E_{21}(r_1)E_{21}(r_2)...E_{21}(r_k).1$ with $k\leq N$.
\endproclaim
\demo{Proof} Prove by induction on the degree of $f$. It is
obvious true for deg $f=0$, i.e $f=1$.

 If $deg f=1$,
$f=x_{(m,n)}=E_{21}(s^m t^n).1$. Now we assume that $f$ is a
monomial whose degree in $x_{(m,n)}$ is positive, then
$$f = x_{(m,n)}f_{\widehat{(m,n)}}=a_{(m,n)}E_{21}(s^m t^n)
f_{\widehat{(m,n)}}-c_{(m,n)}P_{(m,n)}f_{\widehat{(m,n)}}, $$ here
$deg f_{\widehat{(m,n)}}=N-1$. The induction proves our claim.

Hence $E_{21}(r_1)E_{21}(r_2)...E_{21}(r_k).1, k\in
\br_+\bigcup\{0\}, r_i=s^{m_i}t^{n_i},i=1...n,m_i,n_i\in \bz$
spans $V$ over $\bc$. Note that the elements
$E_{21}(r_1)E_{21}(r_2)...E_{21}(r_k).1$ are independent of the
order in which the operators are applied.

Since the leading term of $E_{21}(r_1)E_{21}(r_2)...E_{21}(r_k).1$
is $\prod\Sb i=1\endSb\Sp k\endSp x_{(m_i,n_i)}$, we know that
$$E_{21}(r_1)E_{21}(r_2)...E_{21}(r_k).1$$ form a base for $V$ with
 $k$ ranges in $\{0, 1, 2, 3, \cdots, \}$ and  $r_i$ ranges in $\{s^{m}t^{n}: m, n\in \Bbb{Z} \}$. \qed
\enddemo

It immediately follows from Lemma 4.2 that $V$ is generated as a
$\widetilde{\frak{gl}_{2}(\bc_q)}$-module by $1$, and
$$\bigl(\matrix
  a_1 & a_2 \\
  0 & a_3 \\
\endmatrix\bigr).1=-\frac{1}2\mu \kappa(a_1).1+\frac{1}2\mu
\kappa(a_3).1\tag 4.3$$
 for any $a_1,a_2,a_3\in \bc_q$, here $\kappa(a)$ is defined as in
 $(1.4)$. The subalgebra
 $$\Cal{B}=\{\bigl(\matrix
  a_1 & a_2 \\
  0 & a_3 \\
\endmatrix\bigr): a_1, a_2, a_3 \in \bc_q\} \oplus \bc c_s
\oplus \bc c_t\oplus \bc d_s \oplus \bc d_t$$ is a Borel
subalgebra of $\widetilde{\frak{gl}_{2}(\bc_q)}$ in the sense of
(0.1).

 Hence we have

\proclaim{Proposition 4.4} $V$ is a highest weight module of
highest weight $\lambda: \Cal{B}\to \Bbb{C}$, where $\lambda$ is
defined as follows.
$$\align &\lambda\bigl(\matrix
  a_1 & a_2 \\
  0 & a_3 \\
\endmatrix\bigr)=-\frac{1}2\mu
\kappa(a_1)+\frac{1}2\mu\kappa(a_3), \\
&\lambda(c_s) = \lambda(c_t) =\lambda(d_s) =\lambda(d_t) =
0\endalign$$ and $1$ is the highest weight vector.
\endproclaim

\medskip

Let $i\in \Bbb N$, $\gamma=(\gamma_1,..., \gamma_s)$ be the
$s$-partition  of $i$. Denote by $Par_s(i)$ the set of all
partitions $\gamma=(\gamma_1,..., \gamma_s)$ of $i$ with $s$
parts.

Given $\gamma\in Par_s(N)$, we say that $\pi_1'\times \pi_2'\in
S_N\times S_N$ is equivalent to $\pi_1\times \pi_2\in S_N\times
S_N$, here $S_N$ is the permutation group of $N$ letters, if for
all $z_1,...z_N,w_1,...,w_N\in \bc_q$,
$$\kappa(z_{\pi_1'(1)}w_{\pi_2'(1)}...z_{\pi_1'(\gamma_1)}w_{\pi_2'(\gamma_1)})
...\kappa(z_{\pi_1'(\gamma_1+...\gamma_{s-1}+1)}w_{\pi_2'(\gamma_1+...+\gamma_{s-1}+1)}...
z_{\pi_1'(N)}w_{\pi_2'(N)})$$ can be obtained from the analogous
expression for $\pi_1\times \pi_2$  by a permutation of the $s$
factors $\kappa(\cdots)$ and/or by cyclic permutation of the
variables(e.g.
$\kappa(z_1w_1z_2w_2z_3w_3)=\kappa(z_3w_3z_1w_1z_2w_2)$).

The set of equivalence classes is denoted by $[S_N\times
S_N](\gamma)$. The following result was due to Jakobsen-Kac [JK2].

 \proclaim{Lemma 4.5} Let $z_1,z_2,..z_N,w_1,w_2,..w_N \in
\bc_q[s^{\pm 1},t^{\pm 1}]$
$$\align &\bigl(\matrix
 0 & z_1 \\
  0 & 0 \\
\endmatrix\bigr)\bigl(\matrix
 0 & z_2 \\
  0 & 0 \\
\endmatrix\bigr)...\bigl(\matrix
 0 & z_N \\
  0 & 0 \\
\endmatrix\bigr)\bigl(\matrix
 0 & 0 \\
 w_1  & 0 \\
\endmatrix\bigr)\bigl(\matrix
 0 & 0 \\
 w_2  & 0 \\
\endmatrix\bigr)...\bigl(\matrix
 0 & 0 \\
 w_N  & 0 \\
\endmatrix\bigr).1 \tag 4.6\\
=&\sum\Sb s=1\endSb\Sp N\endSp\sum\Sb \gamma\in
Par_s(N)\endSb\sum\Sb [\pi_1\times\pi_2]\in (S_N\times
S_N)(\gamma)\endSb(-1)^{\gamma_1-1}(-\mu)\kappa(z_{\pi_1(1)}w_{\pi_2(1)}...z_{\pi_1(\gamma_1)}w_{\pi_2(\gamma_1)})\\
&.(-1)^{\gamma_2-1}(-\mu)\kappa(z_{\pi_1(\gamma_1+1)}w_{\pi_2(\gamma_1+1)}...z_{\pi_1(\gamma_2)}w_{\pi_2(\gamma_2)}).\\
&...(-1)^{\gamma_s-1}(-\mu)\kappa(z_{\pi_1(\gamma_1+...\gamma_{s-1}+1)}w_{\pi_2(\gamma_1+...+\gamma_{s-1}+1)}...
z_{\pi_1(N)}w_{\pi_2(N)}).1
\endalign $$
\endproclaim

\medskip

We shall call $k$ the level of the element $E_{21}(r_1)\cdots
E_{21}(r_k).1 \in V$, where $k\in \bz_+\bigcup\{0\}$,
$r_i=s^{m_i}t^{n_i},i=1...n,m_i,n_i\in \bz$.

 \proclaim{Proposition 4.7} (i) The hermitian form on different level is
0.

(ii) Let $h$ be an element of level $n$.  Then $(h, h)$ is a
polynomial in $\mu$ with the leading term $c(h)\mu^n$ for some
constant $c(h)>0$.
\endproclaim
\demo{Proof} Since
$$\align &(1,\bigl(\matrix
 0 & 0 \\
  z_1 & 0 \\
\endmatrix\bigr)\bigl(\matrix
 0 & 0 \\
  z_2 & 0 \\
\endmatrix\bigr)...\bigl(\matrix
 0 & 0 \\
  z_r & 0 \\
\endmatrix\bigr).1)  \ \ \ r\geqq 1\\
=&(\bigl(\matrix
 0 & \overline{z_1} \\
  0 & 0 \\
\endmatrix\bigr).1,\bigl(\matrix
 0 & 0 \\
  z_2 & 0 \\
\endmatrix\bigr)...\bigl(\matrix
 0 & 0 \\
  z_r & 0 \\
\endmatrix\bigr).1)=0
\endalign$$

then WLOG, assume $s>t$,
$$\align &(\bigl(\matrix
 0 & 0 \\
  z_1 & 0 \\
\endmatrix\bigr)\bigl(\matrix
 0 & 0 \\
  z_2 & 0 \\
\endmatrix\bigr)...\bigl(\matrix
 0 & 0 \\
  z_s & 0 \\
\endmatrix\bigr).1,\bigl(\matrix
 0 & 0 \\
  r_1 & 0 \\
\endmatrix\bigr)\bigl(\matrix
 0 & 0 \\
  r_2 & 0 \\
\endmatrix\bigr)...\bigl(\matrix
 0 & 0 \\
  r_t & 0 \\
\endmatrix\bigr).1)\\
=&(1,(-1)^s\bigl(\matrix
 0 & \overline{z_s} \\
  0 & 0 \\
\endmatrix\bigr)\bigl(\matrix
 0 & \overline{z_{s-1}} \\
  0 & 0 \\
\endmatrix\bigr)...\bigl(\matrix
 0 & \overline{z_1} \\
  0 & 0 \\
\endmatrix\bigr)\bigl(\matrix
 0 & 0 \\
  r_1 & 0 \\
\endmatrix\bigr)\bigl(\matrix
 0 & 0 \\
  r_2 & 0 \\
\endmatrix\bigr)...\bigl(\matrix
 0 & 0 \\
  r_t & 0 \\
\endmatrix\bigr).1)\\
=&(1,(-1)^s\bigl(\matrix
 0 & \overline{z_s} \\
  0 & 0 \\
\endmatrix\bigr)\bigl(\matrix
 0 & \overline{z_{s-1}} \\
  0 & 0 \\
\endmatrix\bigr)...\bigl(\matrix
 0 & \overline{z_{s-r}} \\
  0 & 0 \\
\endmatrix\bigr)c.1)\ \ (\text{ according to Lemma 4.5}\text{ here  } c\in
\bc)\\
=&0.
\endalign$$
This proves (i).

\medskip

Let
$h=E_{21}(\overline{z_N})E_{21}(\overline{z_{N-1}})...E_{21}(\overline{z_1})$,
$h'=E_{21}(w_1)E_{21}(w_2)...E_{21}(w_N)$, where
$\overline{z_i}=s^{m_i}t^{n_i}, w_i=s^{l_i}t^{r_i}$, then
according to Lemma 4.5,
$$\align &\overline{(h,h')}\\
=&(-1)^N. \sum\Sb s=1\endSb\Sp N\endSp\sum\Sb \gamma\in
Par_s(N)\endSb\sum\Sb [\pi_1\times\pi_2]\in (S_N\times
S_N)(\gamma)\endSb(-1)^{\gamma_1-1}(-\mu)\kappa(z_{\pi_1(1)}w_{\pi_2(1)}...z_{\pi_1(\gamma_1)}w_{\pi_2(\gamma_1)})\\
&.(-1)^{\gamma_2-1}(-\mu)\kappa(z_{\pi_1(\gamma_1+1)}w_{\pi_2(\gamma_1+1)}...z_{\pi_1(\gamma_2)}w_{\pi_2(\gamma_2)}).\\
&...(-1)^{\gamma_s-1}(-\mu)\kappa(z_{\pi_1(\gamma_1+...\gamma_{s-1}+1)}w_{\pi_2(\gamma_1+...+\gamma_{s-1}+1)}...
z_{\pi_1(N)}w_{\pi_2(N)})\endalign$$ is a polynomial $P$ of $\mu$,
whose coefficients depends on $h$ and $h'$.

If $deg P=N$, then there  exists at least a $\pi\in S_N$, such
that $\kappa(z_iw_{\pi(i)})\neq 0$, that is
$\kappa(t^{-n_i}s^{-m_i}s^{l_{\pi(i)}}t^{r_{\pi(i)}})
=q^{(l_{\pi(i)}-m_i)(-n_i)}\delta_{(l_{\pi(i)}-m_i,r_{\pi(i)}-n_i),(0,0)}\neq
0$, hence $\overline {z_i}=w_{\pi(i)}$for any $1\leq i \leq N$. So
if $h=h'$, the coefficient of $\mu^n$ is the number of such
elements $\pi$, otherwise it equals zero. Hence with Lemma 4.2, we
proved (ii). \qed\enddemo

Next we prove the unitarity of the hermitian form.

\proclaim{Theorem 4.8}$(\pi_{X, \mu},V)$ is unitarizable if and
only if $\mu>0$
\endproclaim
\demo{Proof} Since $$(\bigl(\matrix
 0 & 0 \\
 s^mt^n & 0 \\
\endmatrix\bigr).1,\bigl(\matrix
 0 & 0 \\
  s^lt^k & 0 \\
\endmatrix\bigr).1)=\mu\delta_{m-l,0}\delta_{n-k,0},$$
for any $m,n\in N$, then if $(\pi_{X, \mu},V)$ is unitarizable,
$\mu>0$.

Let $w_i=s^{m_i}t^{n_i},z_j=t^{-l_j}s^{-k_j}$, for $i,j=1,...,N$,
 then $$\align &\kappa(z_1w_1z_2w_2...z_rw_r)\\
 =&\kappa(t^{-l_1}s^{-k_1}s^{m_1}t^{n_1}t^{-l_2}s^{-k_2}s^{m_2}t^{n_2}...t^{-l_r}s^{-k_r}s^{m_r}t^{n_r})\\
 =&q^\alpha \delta_{-k_1+m_1-k_2+m_2-....-k_r+m_r,0}\delta_{-l_1+n_1-l_2+n_2-...-l_r+n_r,0}
 \endalign$$
 where
 $\alpha=(-k_1+m_1)(-l_1)+(-k_2+m_2)(-l_1+n_1-l_2)+...+(-k_r+m_r)(-l_1+n_1-l_2+n_2-...-l_{r-1}+n_{r-1}-l_r)$.
 Consider the linear transformation $T_{a,b}$ of $\bc_q$ determined by
 $$T_{a,b}(s^{c_1}t^{d_1}s^{c_2}t^{d_2}...s^{c_k}t^{d_k})=s^{c_1+a}t^{d_1+b}s^{c_2+a}t^{d_2+b}...s^{c_k+a}t^{d_k+b}$$
 $(a,b\in \bz)$. Extend this operator to a linear
 operator$\widetilde{T_{a,b}}$ on $V$ by
 $$\widetilde{T_{a,b}}[\bigl(\matrix
 0 & 0 \\
  r_1 & 0 \\
\endmatrix\bigr)\bigl(\matrix
 0 & 0 \\
  r_2 & 0 \\
\endmatrix\bigr)...\bigl(\matrix
 0 & 0 \\
  r_t & 0 \\
\endmatrix\bigr).1]=\bigl(\matrix
 0 & 0 \\
  T_{a,b}r_1 & 0 \\
\endmatrix\bigr)\bigl(\matrix
 0 & 0 \\
  T_{a,b}r_2 & 0 \\
\endmatrix\bigr)...\bigl(\matrix
 0 & 0 \\
  T_{a,b}r_t & 0 \\
\endmatrix\bigr).1,$$ for $r_1,...r_t\in \bc_q[s^{\pm 1},t^{\pm 1}]$.

Let $\widetilde{z_i}=T_{-a,-b}(z_i),
\widetilde{w_j}=T_{a,b}(w_j)$, then
$$\align &\kappa(\widetilde{z_1}\widetilde{w_1}\widetilde{z_2}\widetilde{w_2}
...\widetilde{z_r}\widetilde{w_r})\\
 =&\kappa(t^{-l_1-b}s^{-k_1-a}s^{m_1+a}t^{n_1+b}t^{-l_2-b}s^{-k_2-a}s^{m_2+a}t^{n_2+b}...t^{-l_r-b}s^{-k_r-a}
 s^{m_r+a}t^{n_r+b})\\
 =&q^{\widetilde{\alpha}} \delta_{-k_1+m_1-k_2+m_2-....-k_r+m_r,0}\delta_{-l_1+n_1-l_2+n_2-...-l_r+n_r,0}
 \endalign$$
where
 $$\align \widetilde{\alpha}=&(-k_1-a+m_1+a)(-l_1-b)+(-k_2-a+m_2+a)(-l_1-b+n_1+b-l_2-b)+...\\
 &+(-k_r-a+m_r+a)\\
 &(-l_1-b+n_1+b
 -l_2-b+n_2+b-...-l_{r-1}-b+n_{r-1}+b-l_r-b)\\
 =&(-k_1+m_1)(-l_1)+(-k_2+m_2)(-l_1+n_1-l_2)+...\\
 &+(-k_r+m_r)(-l_1+n_1
 -l_2+n_2-...-l_{r-1}+n_{r-1}-l_r)\\
 &-b(-k_1+m_1-k_2+m_2-...-k_r+m_r),
 \endalign$$
 so $\kappa(\widetilde{z_1}\widetilde{w_1}\widetilde{z_2}\widetilde{w_2}
...\widetilde{z_r}\widetilde{w_r})=\kappa(z_1w_1z_2w_2...z_rw_r)$.
It then follows from Lemma 4.2 that $\widetilde{T_{a,b}}$
preserves the hermitian form on $V$.

We need to prove positivity at all levels as the hermitian form on
different levels are zero.

Since $\widetilde{T_{a,b}}$ preserves the hermitian form on $V$,
we may then assume that $h_r$ in level $r$ only involves elements
$s^{m_i}t^{n_i}$ with $m_i\geqq 0,n_i\geqq 0$. Denote
$$\align L_{r}^+(M,N)=Span&\{\bigl(\matrix
 0 & 0 \\
  s^{m_1}t^{n_1} & 0 \\
\endmatrix\bigr)\bigl(\matrix
 0 & 0 \\
  s^{m_2}t^{n_2} & 0 \\
\endmatrix\bigr)...\bigl(\matrix
 0 & 0 \\
  s^{m_r}t^{n_r} & 0 \\
\endmatrix\bigr).1\\
&|m_i\geqq 0, n_i\geqq 0, \sum\Sb i=1\endSb\Sp r\endSp m_i\leqq M,
\sum\Sb i=1\endSb\Sp r\endSp n_i\leqq N\}; \endalign$$

  From the above discussions, we know that the
hermitian form restricted to every level should be positive
definite for $\mu$ big enough. Assume that the form is not
positive definite for some ( possible all ) $\mu > 0$. Let $s_0$
be the lowest level at which there is non-unitarity. It is clear
that $s_0>1$. So there exist $M,N$ such that the form restricted
to $L_{s_0}^+(M,N)$ is not positive definite. Following (4.7), the
form on $L_{s_0}^+(M,N)$ varies smoothly with $\mu$, then we can
find a $\mu_0$ ( we can think it is the first place going from
$\infty$ towards 0 ) at which the form is not positive definite,
while for all $\mu>\mu_0$, the form is positive definite. We write
$(. , .)_{\mu}$ to be the hermitian form at $\mu$.
\medskip
\proclaim {Claim} The radical of the form is non-trivial at
$\mu_0$.
\endproclaim

At first for all $h'\in L_{s_0}^+(M,N)$, $(h',h')_{\mu_0}\geq 0$.
Otherwise, there exists $h\in L_{s_0}^+(M,N)$ such that
$(h,h)_{\mu_0}<0$. From (4.7), the form varies smoothly with
$\mu$, and $(h,h)_{\mu}>0$ if $\mu\rightarrow \infty$, then there
exist $\mu'>\mu_0$ such that $(h,h)_{\mu'}=0$, this contradicts
with the fact that for all $\mu>\mu_0$, the form is positive
definite.

Since the form is positive semi-definite but not positive definite
at $\mu_0$, the radical of the form must be non-trivial.
 Thus,
$$\exists 0\neq\widetilde{h}\in L_{s_0}^+(M,N), \forall h\in L_{s_0}^+(M,N) \text{ such that }(
\widetilde{h}, h)_{\mu_0}=0.$$
 Let $h_{s_0-1}$ be an arbitrary element of
$L_{s_0-1}^+(M,N)$, and let $c\in \bc$, then
$$(\bigl(\matrix
 0 & c \\
  0 & 0 \\
\endmatrix\bigr)\widetilde{h},h_{s_0-1})_{\mu_0}=0.$$

>From the assumption of $s_0$, we have $\bigl(\matrix
 0 & c \\
  0 & 0 \\
\endmatrix\bigr)\widetilde{h}=0$, for any $c\in \bc$.
Replacing $\widetilde{h}$ by
$\widetilde{T_{-m,-n}}(\widetilde{h})$ if necessary, we can write
$$\widetilde{h}=\sum\Sb i=1\endSb\Sp s_0\endSp a_i(\bigl(\matrix
 0 & 0 \\
  1 & 0 \\
\endmatrix\bigr)^i x_i.1$$ where $x_i=\sum \prod_{j=1}^{s_0-i}\bigl(\matrix
 0 & 0 \\
  v_{i,j} & 0 \\
\endmatrix\bigr)$ (here is the finite sum), and each $v_{i,j}$ is the form of $s^lt^k$ (here $l, k$ can not both be $0$).

Let $i_0$ be the smallest $i, 1\leq i \leq s_0$ such that
$a_{i_0}\neq 0$. It follows that
$$\bigl(\matrix
 0 & c \\
  0 & 0 \\
\endmatrix\bigr)\widetilde{h}=\beta a_{i_0}\bigl(\matrix
 0 & 0 \\
  1 & 0 \\
\endmatrix\bigr)^{i_0-1} x_{i_0}.1+R=0$$
where $R$ contains a power of $\bigl(\matrix
 0 & 0 \\
  1 & 0 \\
\endmatrix\bigr)$  greater than $i_0-1$. Observe that
$$\align \bigl(\matrix
 0 & c \\
  0 & 0 \\
\endmatrix\bigr)\bigl(\matrix
 0 & 0 \\
  1 & 0 \\
\endmatrix\bigr)^{i_0}=&\bigl(\matrix
 0 & 0 \\
  1 & 0 \\
\endmatrix\bigr)^{i_0}\bigl(\matrix
 0 & c \\
  0 & 0 \\
\endmatrix\bigr)+i_0\bigl(\matrix
 0 & 0 \\
  1 & 0 \\
\endmatrix\bigr)^{i_0-1}\bigl(\matrix
 c & 0 \\
  0 & -c \\
\endmatrix\bigr)\\
&+(-2c)\frac{i_0(i_0-1)}{2}\bigl(\matrix
 0 & 0 \\
  1 & 0 \\
\endmatrix\bigr)^{i_0-1}\endalign$$
(this can be  easily proved by induction on $i_0$).

Since$$\bigl(\matrix
 c & 0 \\
  0 & -c \\
\endmatrix\bigr)x_{i_0}.1=x_{i_0}\bigl(\matrix
 c & 0 \\
  0 & -c \\
\endmatrix\bigr).1+(-2c)(s_0-i_0)x_{i_0},$$
we have
$$\beta=c(-i_0\mu_0+i_0(-2)(s_0-i_0))+(-c)i_0(i_0-1)=ci_0(-\mu_0-(s_0-i_0)-(s_0-1)).$$

Since $s_0\geqq i_0\geqq 1$ and $\mu_0>0$, $\beta\neq 0$ which
contradicts with $\bigl(\matrix
 0 & c \\
  0 & 0 \\
\endmatrix\bigr)\widetilde{h}=0$.

So for any $\mu > 0$, the hermitian form is positive definite.
\qed\enddemo

\bigskip

\heading{\bf Acknowledgments}\endheading

\medskip

We would like to thank Professor Minoru Wakimoto for sending us
the unpublished manuscript [W1]. We are also grateful to
Professors Bruce Allison and Chongying Dong for some helpful
conversations. A word of special thanks goes to the referee who
provided many suggestions which significantly improved the
exposition of this paper.

\enddocument